\documentclass[journal,english]{IEEEtran}

\usepackage{graphicx}
\usepackage{acronym}
\usepackage{amsmath}
 \usepackage{url,doi}

\usepackage{SIunits}
\usepackage{balance}
\usepackage[linesnumbered]{algorithm2e}

\usepackage[T1]{fontenc}
\usepackage{placeins}
\usepackage{color}
\usepackage[normalem]{ulem}

\newcommand{\state}{x}
\newcommand{\mean}{\mu}

\newcommand{\priormean}{\mean^-}

\newcommand{\msqrt}[1]{\sqrt{ #1}}

\newcommand{\el}[1]{_{[#1]}}

\newcommand{\BSloc}{r}

\newcommand{\SPLgamma}{\gamma}

\newcommand{\statedim}{n}
\newcommand{\measdim}{d}
\newcommand{\measfun}{h}
\newcommand{\meas}{y}
\newcommand{\statenoise}{\varepsilon^\state}
\newcommand{\measnoise}[1][]{\varepsilon^\meas}
\newcommand{\hatmeasnoise}[1][]{\hat{\varepsilon}^\meas}
\newcommand{\meascov}{R}
\newcommand{\hatmeascov}{\hat{R}}

\newcommand{\statecov}{P}
\newcommand{\priorcov}{\statecov^-}

\newcommand{\nonl}{\eta}

\newcommand{\nonllimit}{\nonl_\text{threshold}}
\newcommand{\statefun}{f}
\newcommand{\stateerrcov}{W}

\newcommand{\eigVALUE}{\Lambda}
\newcommand{\eigVECTOR}{U}

\newcommand{\trvec}{\xi^h}
\newcommand{\hattrvec}{\hat{\xi}^h}

\newcommand{\trmat}{\Xi^h}
\newcommand{\hattrmat}{\hat{\Xi}^h}

\newcommand{\trvs}{\xi^f}
\newcommand{\trms}{\Xi^f}

\newcommand{\truestate}{\state_\text{true}}

\newcommand{\decmat}{D}

\newcommand{\innocov}{S}
\newcommand{\kalmangain}{K}

\newcommand{\predmeas}{\meas^-}
\newcommand{\hatpredmeas}{\hat{\meas}^-}

\newcommand{\jacobian}{{J^h}}

\newcommand{\jacobianh}{{J^h}}
\newcommand{\jacobianf}{{J^f}}
\newcommand{\hessian}{H^h}
\newcommand{\hathessian}{\hat{H}^h}

\newcommand{\hessianh}{H^h}

\newcommand{\hessians}{H^f}

\newcommand{\sqvec}{\Delta}

    \providecommand{\matr}[1]{\begin{bmatrix} #1 \end{bmatrix}}

    \DeclareMathOperator{\N}{N}                                 
    \DeclareMathOperator{\tr}{tr}                                 
    \DeclareMathOperator{\atan}{atan2}

\begin{document}
\newacro{AP}{Access Point}
\newacro{APLS}{Access Point Least Squares}
\newacro{AS}{Adaptive Splitting}
\newacro{BinoGM}{Binomial Gaussian Mixture}
\newacro{BinoGMF}{Binomial Gaussian Mixture Filter}
\newacro{BS}{Basestation}
\newacro{CA}{Coverage Area}
\newacro{cdf}{cumulative distribution function}
\newacro{CDF}{Central Difference Filter}
\newacro{CKF}{Cubature Kalman Filter}
\newacro{EKF}{Extended Kalman Filter}
\newacro{EKF2}{Second Order Extended Kalman Filter}
\newacro{IPLF}{Iterated Posterior Linearisation Filter}
\newacro{FP}{Fingerprint}
\newacro{GGM}{Generalized Gaussian Mixture}
\newacro{GM}{Gaussian Mixture}
\newacro{GMEM}{Gaussian Mixture Expectation Maximization}
\newacro{GMF}{Gaussian Mixture Filter}
\newacro{GN}{Gauss-Newton}
\newacro{GNMax}{Gauss-Newton Max  Range}
\newacro{GNSS}{Global Navigation Satellite System}
\newacro{GPS}{Global Positioning System}
\newacro{GSM}{Global System for Mobile Communications}
\newacro{IKF}[IEKF]{Iterated Extended Kalman Filter}
\newacro{IMU}{Inertial Measurement Unit}
\newacro{KL}{Kullback-Leibler}
\newacro{LTE}{Long Term Evolution}
\newacro{OEKF2}[PUKF]{Partitioned Update Kalman Filter}
\newacro{pdf}{probability density function}
\newacro{PF}{Particle Filter}
\newacro{PDR}{Pedestrian Dead Reckoning}
\newacro{PL}{Pathloss}
\newacro{ROTHEKF2}[NEKF2]{Numerical Second Order Extended Kalman Filter}
\newacro{RSS}{Received Signal Strength}
\newacro{RUF}{Recursive Update Filter}
\newacro{TA}{Timing Advance}
\newacro{TDoA}{Time Difference of Arrival}
\newacro{UKF}{Unscented Kalman Filter}
\newacro{VARAPLS}{Variance Access Point Least Squares}
\newacro{WKNN}{Weighted $k$-nearest Neighbor}
\newacro{WLAN}{Wireless Local Area Network}

\title{Partitioned Update Kalman Filter}
\author{  \IEEEauthorblockN{Matti Raitoharju\IEEEauthorrefmark{1}, Robert Piché\IEEEauthorrefmark{1}, Juha Ala-Luhtala\IEEEauthorrefmark{2}, Simo Ali-Löytty\IEEEauthorrefmark{2}} \\ \{matti.raitoharju, robert.piche, juha.ala-luhtala, simo.ali-loytty\}@tut.fi \thanks{
 \IEEEauthorblockA{\IEEEauthorrefmark{1}Department of Automation Science and Engineering, Tampere University of Technology} 
    \IEEEauthorblockA{\IEEEauthorrefmark{2}Department of Mathematics, Tampere University of Technology.
     The authors declare that they have no competing interests.
    Juha Ala-Luhtala received financial support from the Tampere University of Technology Doctoral Programme in Engineering and Natural Sciences.
    The simulations were carried out using the computing resources of CSC --- IT Center for Science.
}
  }  }
\maketitle
\allowdisplaybreaks

\begin{abstract}
In this paper we present a new Kalman filter extension for state update called \ac{OEKF2}. \ac{OEKF2} updates the state using multidimensional measurements in parts. \ac{OEKF2} evaluates the nonlinearity of the measurement function within a Gaussian prior by comparing the effect of the 2nd order term on the Gaussian measurement noise. A linear transformation is applied to measurements to minimize the nonlinearity of a part of the measurement. The measurement update is then applied using only the part of the measurement that has low nonlinearity and the process is then repeated for the updated state using the remaining part of the transformed measurement until the whole measurement has been used. \ac{OEKF2} does the linearizations numerically and no analytical differentiation is required. Results show that when the  measurement geometry allows effective partitioning, the proposed algorithm improves estimation accuracy and produces accurate covariance estimates.
\end{abstract}

\acresetall
\section{Introduction}
Bayesian filtering algorithms are used to compute the estimate of an $\statedim$-dimensional state $\state$. In a general discrete-time model the state evolves according to a  state transition equation
\begin{align}
	\state_{t} = \tilde\statefun_t\left(\state_{t-1},\statenoise_t\right), \label{equ:ssg}
\end{align} 
where $\tilde\statefun_t$ is the state transition function at time index $t$ and $\statenoise_t$ is the state transition noise. The state estimate is updated using measurements that are modeled as
\begin{equation}
	\meas_t=\tilde\measfun_t( \state_t, \measnoise_t) \label{equ:meg},
\end{equation}
where $\tilde\measfun_t$ is a measurement function and $\measnoise_t$ is the measurement noise. If the measurement  and  state transition  are linear,  noises are additive, white and normal distributed, and the prior state ($\state_0$) is normal distributed, the Kalman update can be used to compute the posterior. If these requirements are not fulfilled, usually an approximate estimation method has to be used. In this work, we concentrate on situations where the noises are additive and Gaussian so that (\ref{equ:ssg}-\ref{equ:meg}) take the form
\begin{align}
	\state_{t} &= \statefun_t\left(\state_{t-1}\right) + \statenoise_t \label{equ:addtrans}\\
	\meas_t &=\measfun_t( \state_t ) + \measnoise_t, \label{equ:measurement}
\end{align}
where  $\statenoise_t \sim \N(0,\stateerrcov_t)$, $\stateerrcov?_t$ is the state transition noise covariance,  $\measnoise_t \sim \N(0,\meascov_t)$, and $\meascov_t$ is the measurement noise covariance.

There are two main approaches for computing an approximation of the posterior distribution:
\begin{enumerate}
	\item Approximate probabilities using point masses (e.g.\ grid and particle filters)
	\item Approximate probabilities by Gaussians (e.g.\ Kalman filter extensions)
\end{enumerate}
In the first approach one problem is how to choose a good number of point masses. The first approach also often requires more computational resources than the second approach. A drawback of the second approach is that the state distribution is assumed normal and unimodal, which makes the estimate inaccurate when the true posterior is not normal. \acp{GMF} (a.k.a. Gaussian sum filters) can be considered as a hybrid approach that use multiple normal distributions to estimate the probability distributions and can approximate any \ac{pdf}~\cite{Sorenson_Alspach_1971}. \acp{GMF}\ have the same kind of problems as the algorithms using point masses in choosing a good number of components. The algorithm that will be proposed in this paper uses the second approach and so we will concentrate on it.

The algorithms that are based on Gaussian approximations usually extend the Kalman filter update to nonlinear measurements (there are also other options, see for example \cite{gf1,gf2}).
The \ac{EKF} is a commonly used algorithm for estimation with nonlinear measurement models \cite{gelb1974aoe}. \ac{EKF} is based on the first order Taylor linearization of the measurement function at the mean of prior. In the \ac{EKF2} the linearization takes also the second order expansion terms into account \cite{gelb1974aoe}. In contrast to \ac{EKF}, in \ac{EKF2} the prior covariance also affects the linearization. Both \ac{EKF} and \ac{EKF2} require analytical computation of the Jacobian matrix and \ac{EKF2} requires also the computation of Hessian matrices of the measurement function.  In \cite{855552} a 2nd order \ac{CDF}, which can be interpreted as a derivative-free numerical approximation of \ac{EKF2}, was presented. The most commonly used Kalman filter extension that does not require analytical differentiation is probably the \ac{UKF} \cite{julier1997new}. The Gaussian approximations in \ac{UKF} are based on the propagation of ``sigma points'' through the nonlinear functions. \acp{CKF} are similar algorithms, but they have different theory in the background \cite{cubature}. All these methods do the update as a single operation.

Some algorithms do multiple linearizations to improve the estimate. In~\cite{dmitriev1978generalized} the posterior is computed using multiple \ac{EKF} updates that use different linearization points. In the \ac{IKF}, the \ac{EKF} update is computed in the prior mean and then the new mean is used as the new linearization point \cite{jazwinski}. This can be done several times. A similar update can be done also with other Kalman type filters \cite{6916133}. The \ac{RUF} updates the prior with measurement with reduced weight several times~\cite{RUKF}. In every update the linearization point is used from the posterior of the last reduced weight update. \acp{GMF} can also be considered to be filters that do the linearization multiple times, once for each Gaussian component, and any Kalman filter extension can be used for the update.
 
In this paper we present \ac{OEKF2} that updates the state also in several steps. \ac{OEKF2} first computes the nonlinearity of  measurement models. The nonlinearity measure is based on comparing the covariance of the 2nd order term covariance of the Gaussian measurement noise. Computation of this nonlinearity measure requires the same matrices as the \ac{EKF2} update and for this we use the 2nd order \ac{CDF} \cite{855552}, which is a derivative free version of the \ac{EKF2}.

\begin{figure}
	\includegraphics[width=\columnwidth]{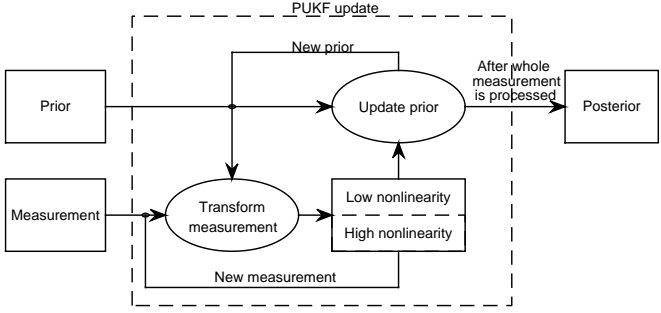}
	\caption{Process diagram of the PUKF}
	\label{fig:diagram}
\end{figure}

\ac{OEKF2} applies a linear transformation to the measurement function to make a new measurement function that has linearly independent measurement noise for measurement elements; the smallest nonlinearity corresponding to a measurement element is minimized first, then the second smallest nonlinearity etc.  After the transformation, the update is done using only measurement elements that have smaller nonlinearity than a set threshold value or using the measurement element with the smallest nonlinearity. After the partial measurement update the covariance has become smaller or remained the same and the linearization errors for remaining measurements may have also became smaller. The remaining measurements' nonlinearity is re-evaluated using the partially updated state, the remaining measurements are transformed and a new partial update is applied until the whole measurement is applied. This process is shown in Figure~\ref{fig:diagram}. The use of only some dimensions of the measurements to get a new prior and the optimization of measurement nonlinearities differentiates \ac{OEKF2} from other Kalman filter extensions.

The article is structured as follows: In Section~\ref{sec:numup} a numerical method for approximate \ac{EKF2} update is presented. The main algorithm is presented in Section~\ref{sec:PEKF}. The accuracy and reliability of the proposed algorithm is compared with other Kalman filter extensions and \acp{PF} in Section~\ref{sec:tests}. Section~\ref{sec:conclusions} concludes the article.

\section{\ac{EKF2} and Its Numerical Update Using 2nd Order \ac{CDF}}
\label{sec:numup}
Kalman filter extensions, like all Bayesian filters, can be computed in two stages: prediction and update. For the state transition model (\ref{equ:addtrans}) the state is propagated  in \ac{EKF2} using equations \cite{jazwinski}:
\begin{align}
	\priormean_t & = \statefun_t(\mean^+_{t-1}) + \frac{1}{2} \trvs_t \label{equ:pred1} \\
	\priorcov_t &= \jacobianf \statecov^+_{t-1}  \jacobianf^T + \frac{1}{2}\trms_t + \stateerrcov_t, \label{equ:pred2} 
\end{align}
where $\priormean_t$ is the predicted mean at time $t$, $\mean^+_{t-1}$ is the posterior mean of the previous time step, $\jacobianf$ is the Jacobian of the state transition function evaluated at $\mean^+_{t-1}$, $\priorcov_t$ is the predicted covariance, $\statecov^+_{t-1}$ is the posterior covariance of the previous time step and $\trvec_t$ and $\trmat_t$ are defined as
\begin{align}
	{\trvs_t}\el{i} & =  \tr \statecov^+_{t-1}\hessians_i \label{equ:trsvec} \\
	{\trms_t}\el{i,j} & =  \tr \statecov^+_{t-1}\hessians_i\statecov^+_{t-1}\hessians_j, \label{equ:trsmat}
\end{align}
where $\hessians_i$ is the Hessian of the $i$th element of the state transition function evaluated at $\mean^+_{t-1}$. 
To simplify the notation we do not further show the time indices.

The update equations of \ac{EKF2} for the measurement model (\ref{equ:measurement}) are \cite{jazwinski}
\begin{align}
	\predmeas & = \measfun(\priormean)+ \frac{1}{2} \trvec \label{equ:predmeas} \\
	\innocov & = \jacobianh \priorcov \jacobianh^T + \frac{1}{2} \trmat +  \meascov \label{equ:innocov}	\\
	\kalmangain & = \priorcov\jacobianh^T \innocov^{-1} \label{equ:Kgain}\\
	\mean^+ &= \priormean + \kalmangain ( \meas - \predmeas) \\
	\statecov^+ &=  \priorcov - \kalmangain\innocov\kalmangain^{T}, \label{equ:ekfend}
\end{align}
where $\jacobianh$ is the Jacobian of the measurement function, $\kalmangain$ is the Kalman gain, $\innocov$ is the innovation covariance, and $\trvec$ and $\trmat$ are defined as
\begin{align}
	\trvec\el{i} & =  \tr \priorcov\hessianh_i \label{equ:trvec} \\
	\trmat\el{i,j} & =  \tr \priorcov\hessianh_i\priorcov\hessianh_j, \label{equ:trmat}
\end{align}
where $\hessianh_i$ is the Hessian matrix of the $i$th component of the measurement function. Eqns (\ref{equ:predmeas}-\ref{equ:ekfend}) can be turned into the \ac{EKF} update using $\trvec=0$ and $\trmat=0$.

If the measurement model is linear, the trace terms in \ac{EKF2} are zero and the update is the optimal update of the Kalman filter. When the measurement function is a second order polynomial the \ac{EKF2} update is not optimal as the distributions are no longer Gaussian, but the mean (\ref{equ:predmeas}) and innovation covariance (\ref{equ:innocov}) are correct. 

In this paper we use a numerical algorithm to compute an \ac{EKF2} like update. To derive this algorithm, we start with the formulas of the 2nd-order CDF from \cite{855552}. 
 Let $\msqrt{\priorcov}$ be a matrix such that
\begin{equation}
	\msqrt{\priorcov}\msqrt{\priorcov}^T=\priorcov.  \label{equ:sqrtp}
\end{equation}
In our implementation this matrix square root is computed using Cholesky decomposition. 

Next we define matrices $M$ and $Q$ that are used for computing the numerical \ac{EKF2} update. We use notation $\sqvec_i=\SPLgamma\msqrt{\priorcov}\el{:,i}$, where $\msqrt{\priorcov}\el{:,i}$ is the $i$th column of matrix $\msqrt{\priorcov}$ and $\SPLgamma$ is an algorithm parameter that defines the spread of the function evaluations. Matrix $M$, whose elements are
\begin{equation}
\begin{aligned}
		M\el{:,i}=&\left[\jacobian\msqrt{\priorcov}\right]\el{:,i} \\ 
		\approx &\SPLgamma^{-1}\frac{\measfun(\priormean +\sqvec_i ) -  \measfun(\priormean - \sqvec_i )}{2},
		\end{aligned}
	\label{equ:M}
\end{equation} is needed for the terms with Jacobian.
The matrices $Q_k\approx\msqrt{\priorcov}\hessian_k\msqrt{\priorcov}^T$ are needed to compute terms with Hessians. Elements of $Q_k$ are
	\begin{align}
	\begin{aligned}
		{Q_k}\el{i,i} = &\SPLgamma^{-2}\left[\measfun\el{k}(\priormean +\sqvec_i ) +  \measfun\el{k}(\priormean - \sqvec_i) - 2\measfun\el{k}(\priormean )\right] \\
		{Q_k}\el{i,j} = &\SPLgamma^{-2}\left[\measfun\el{k}(\priormean +\sqvec_i +\sqvec_j) - \measfun\el{k}(\priormean +\sqvec_i )\right. \\ 
	&	\left. -\measfun\el{k}(\priormean +\sqvec_j ) + \measfun\el{k}(\priormean)\right], i\neq j. \label{equ:Q}
	\end{aligned}
	\end{align}
The \ac{EKF2} update can be approximated with these by doing the following substitutions:
	\begin{align}
		\trvec_i=\tr \priorcov\hessian_i & \approx \tr Q_i && \text{in (\ref{equ:predmeas})} \label{equ:Q1} \\
		\jacobian\priorcov\jacobian^T&  \approx  MM^T && \text{in (\ref{equ:innocov})}\\
		\priorcov\jacobian^T & \approx \msqrt{\priorcov} M^T && \text{in (\ref{equ:trvec})}\\
		\trmat\el{i,j}=\tr \priorcov\hessian_i\priorcov\hessian_j & \approx \tr Q_iQ_j \label{equ:QQ} && \text{in (\ref{equ:trmat})}.
	\end{align}
	The prediction step can be approximated by computing $M^f$~(\ref{equ:M}) and $Q^f$~(\ref{equ:Q}) matrices  using the state transition function instead of the measurement function and doing the following substitutions:
	\begin{align}
		\jacobian\statecov_{t-1}^+\jacobian^T&  \approx  M^f{M^f}^T && \text{in (\ref{equ:pred2})}\\
		\tr \priorcov\hessians_i & \approx \tr Q^f_i && \text{in (\ref{equ:trsvec})} \\
		\tr \priorcov\hessians_i\priorcov\hessians_j & \approx \tr Q^f_iQ^f_j && \text{in (\ref{equ:trsmat})}.
	\end{align}
	
In \cite{norgaard}, an update algorithm similar to numerical \ac{EKF2} is proposed that uses only the diagonal elements of $Q$ matrices. They state that $\SPLgamma=\sqrt{3}$ for Gaussian distributions is optimal because it preserves the fourth moment and so we use this $\SPLgamma$ value in our algorithm.
	
\section{Partitioned Update Kalman Filter}
\label{sec:PEKF}
When the measurement function is linear and the measurement noise covariance is block diagonal, the Kalman update produces identical results whether measurements are applied one block at a time or all  at once. In our approach we try to find as linear as possible part of the measurement and use this part to update the state estimate to reduce approximation errors in the remaining measurement updates. When the measurement noise covariance $\meascov$ is  not diagonal a linear transformation (decorrelation) is applied to transform the measurement so that the transformed measurement has diagonal covariance \cite{stengel1986stochastic}. In \ac{OEKF2}, we choose this decorrelation so that the nonlinearity of the least nonlinear measurement element is minimized. The prior is updated using the least nonlinear part of the decorrelated measurements. After the partial update the process is repeated for the remaining dimensions of the transformed measurement.

For measuring the amount of nonlinearity we compare the trace term $\trmat$ with the covariance of the measurement noise:
\begin{equation}
\begin{aligned}
	\nonl&=\tr \sum_{i=1}^\measdim \sum_{k=1}^\measdim \meascov\el{k,i}^{-1} \priorcov \hessian_i  \priorcov \hessian_k  \label{equ:nonl} \\
	&= \tr  \sum_{i=1}^\measdim \sum_{k=1}^\measdim \meascov\el{k,i}^{-1} \trmat\el{i,k}
\end{aligned}
\end{equation} 
This nonlinearity measure is a local approximation of the nonlinearity and is developed from the measure presented in \cite{jazwinski,ali-loytty2010c}. In~\cite{raitoharjuphd} it was compared with other nonlinearity measures and it was shown to be a good indication of how accurately state can be updated with a nonlinear measurement model using a Kalman filter extension. When the measurement model is linear the nonlinearity measure is $\nonl=0$.

\begin{figure}
	\includegraphics[width=\columnwidth,clip=true,trim=1cm 0cm 1cm 0cm]{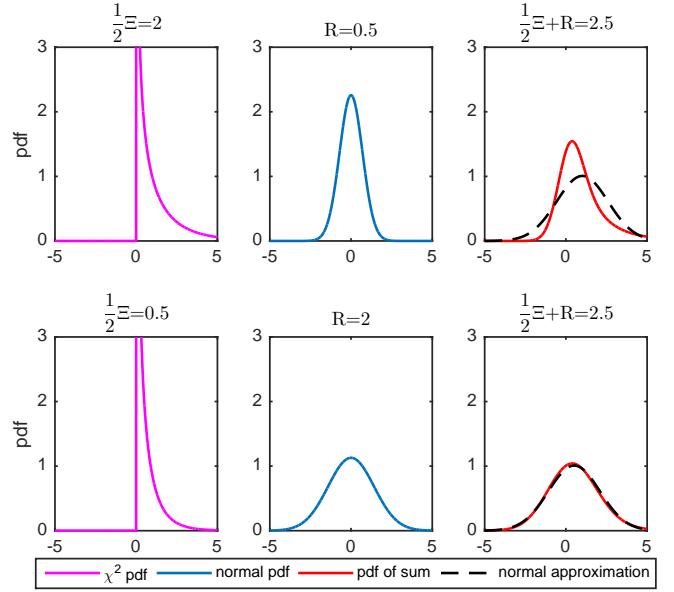}
	\caption{Probability density functions of sums of independent $\chi^2$ and normal random variables with different variances}
	\label{fig:nonlinex}
\end{figure}

The matrix $\trmat$ depends on the nonlinearity of the measurement function and contributes to the innovation covariance (\ref{equ:innocov}) similarly, but multiplied with $\frac{1}{2}$, as the Gaussian measurement noise $R$. The measure (\ref{equ:nonl}) compares the ratio of Gaussian covariance $R$ and non-Gaussian covariance $\trmat$.  Figure~\ref{fig:nonlinex} shows how the pdf of the sum of independent normal and $\chi^2$ distributed random variables is closer to normal when $R>\Xi$ than when $R<\Xi$. The $\chi^2$ distribution is chosen in the example, because a normal distributed variable squared is $\chi^2$  distributed and in the second order polynomial approximations the squared term is the nonlinear part.

 The nonlinearity measure (\ref{equ:nonl}) can be approximated numerically using the substitution (\ref{equ:QQ}). Numerical computation of a similar nonlinearity measure was proposed in \cite{AS}, but the algorithm presented in Section~\ref{sec:numup} does the nonlinearity computation with  fewer measurement function evaluations. 

Multiplying (\ref{equ:measurement}) by an invertible square matrix $\decmat$ gives a transformed measurement model
\begin{align}
	\decmat \meas & = \decmat \measfun(\state)+\decmat \measnoise.
\end{align}
We use the following notations for the transformed measurement model: $\hat{\meas}=\decmat \meas$, $\hat{\measfun}(x)=\decmat {\measfun}(x)$, $\hat{R}=\decmat R\decmat^T$, and ${\hatmeasnoise=D\measnoise \sim \N(0,\hat R)}$. We will show that  $\decmat$ can be chosen so that 
\begin{equation}
	 \hat{\meascov} = I \quad \text{and}\quad \tr{\priorcov\hathessian_i\priorcov\hathessian_k} = 0 , i \neq k ,
\end{equation}
where $\hathessian_i$ and $\hathessian_k$ denote the Hessians of the $i$th and $k$th element of $\hat{h}(x)$. 

In~\cite{raitoharjuphd}, it was shown that when a measurement model is transformed so that $\hat{R}=I$ the nonlinearity measure (\ref{equ:nonl}) is equal to the nonlinearity measure of the transformed measurements
 \begin{equation}
 \begin{aligned}
\eta&=\tr \sum_{i=1}^\measdim \sum_{k=1}^\measdim \meascov\el{k,i}^{-1} \priorcov \hessian_i  \priorcov \hessian_k   \\
&=\tr \sum_{i=1}^\measdim \sum_{k=1}^\measdim \hat\meascov\el{k,i}^{-1} \priorcov \hathessian_i  \priorcov \hathessian_k  =\tr \sum_{i=1}^\measdim \priorcov \hathessian_i  \priorcov \hathessian_i.
\end{aligned}
\end{equation} 
Because the cross terms do not affect to the amount of nonlinearity we can extract the nonlinearity caused by individual elements of the transformed measurements
\begin{equation}
\nonl_i=\tr \priorcov \hathessian_i  \priorcov \hathessian_i \label{equ:ithnonl}
\end{equation}
and the total nonlinearity is 
\begin{equation}
\nonl=\sum_{i=1}^\measdim \nonl_i .
\end{equation}

In Appendix~\ref{app:lintrans} it is shown that 
\begin{equation}
	\hattrmat\el{i,j} = \left[\decmat \trmat \decmat^T\right]\el{i,j} \approx  \tr{\priorcov\hathessian_i\priorcov\hathessian_j}.
\end{equation} In this case the measurement related error terms of the transformed measurement $\hatmeascov$ and $\hattrmat$ are diagonal. This makes the measurements independent and allows the update of the state one element at a time.

In \ac{OEKF2}  nonlinearities are minimized in such a way that $\eta_1$ (\ref{equ:ithnonl}) is as small as possible. Then $\eta_2$ is minimized such that $\eta_1$ does not change, and $\eta_3$ so that $\eta_1$ and $\eta_2$ do not change etc. The decorrelation transformation $\decmat$ that does the desired nonlinearity minimization can be computed by first computing a matrix square root (\ref{equ:sqrtp}) of the measurement noise covariance   
\begin{equation}
	\msqrt{\meascov}\msqrt{\meascov}^T = \meascov \label{equ:sqrtr}
\end{equation}
and then an eigendecomposition of $\msqrt{\meascov}^{-1} \trmat \msqrt{\meascov}^{-T}$
\begin{equation}
	\eigVECTOR\eigVALUE\eigVECTOR^T = \msqrt{\meascov}^{-1} \trmat \msqrt{\meascov}^{-T}. \label{equ:eigdecom}
\end{equation}
We assume that the eigenvalues in the diagonal matrix $\eigVALUE$ are sorted in ascending order. 
The transformation matrix is now
\begin{equation}
	\decmat = \eigVECTOR^{T} \msqrt{\meascov}^{-1}. \label{equ:D}
\end{equation}
A proof that this transformation minimizes the nonlinearity measures is given in Appendix~\ref{app:nonlminimization}. After transforming the measurement model with this matrix, the measurement noise covariance is $\hatmeascov=I$ and $\hattrmat=\eigVALUE$. 

After the measurement model is decorrelated (multiplied with~$\decmat$), the parts of measurement model that have low nonlinearity ($\Lambda\el{i,i} \leq \nonllimit$) are  used in the update (Section~\ref{sec:numup}). If there is no such part then the most linear element of the measurement model is used to update the state. Then the same process is repeated for the remaining transformed measurement model until the whole measurement is processed. 

In summary the \ac{OEKF2} update is:
\begin{enumerate}
	\item Transform the measurement model using $D$ (\ref{equ:D})
	\item Update the prior using only the least nonlinear measurement elements of the transformed measurement
	\item If there are measurement elements left, use them as new measurement and use the updated state as a new prior and return to step 1
\end{enumerate}
The detailed \ac{OEKF2} algorithm is presented in Algorithm~\ref{algo:OEKF} and a Matlab implementation is available online \cite{matlab}.

The amount of nonlinearity (\ref{equ:nonl}) for independent measurements is equal to the sum of the nonlinearities for each of the measurements. The quantity $\nonllimit$ is compared separately to independent transformed measurements elements and, thus, we propose to use same $\nonllimit$ regardless of the measurement dimension. As a rule of thumb the nonlinearity threshold can be set to $\nonllimit=1$, which is equal to the threshold proposed for one dimensional measurements in \cite{jazwinski}.

\begin{figure}[!tpb]
	\includegraphics[width=\columnwidth,clip=true,trim=1.2cm 0.5cm 1.2cm 0.5cm]{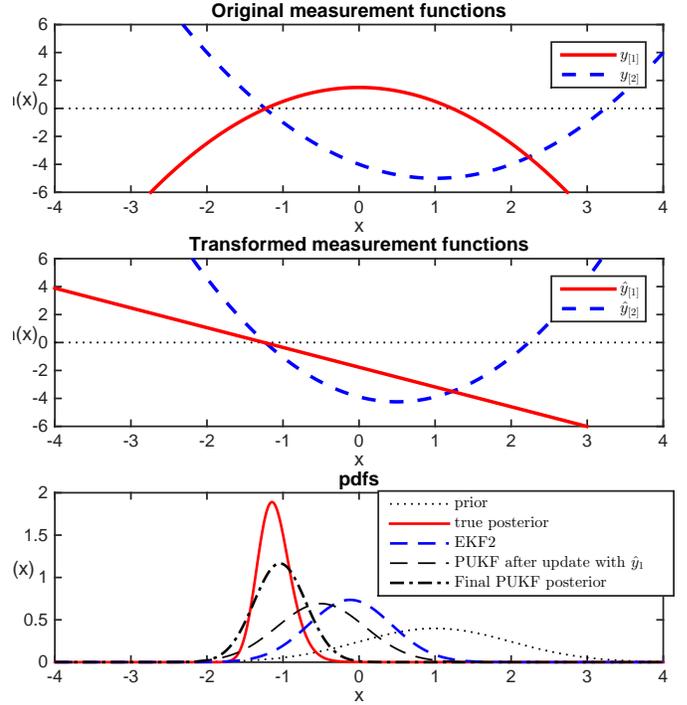}
	\caption{Transforming second order polynomial measurements to minimize nonlinearity of $\hat{\meas}_1$ and posterior comparison of  \ac{OEKF2} and \ac{EKF2}}
	\label{fig:pdfs}
\end{figure}

\begin{algorithm}[tb]
       \SetKwInOut{Input}{input}
        \SetKwInOut{Output}{output}
        \SetKwInOut{Parameters}{parameters}
        \Input{
        Prior state: $\mean$ --  mean 
                  $\statecov$ --  covariance \newline
        Measurement model: $\meas$ --  value, $\measfun(\cdot)$ --  function, 
                  $\meascov$~--~covariance \newline
                   $\nonllimit$ -- nonlinearity limit, $\SPLgamma$~--~measurement function evaluation spread (default $\SPLgamma=\sqrt{3}$)
                  }
         \Output{
         Updated state: $\mean$ -- mean, $\statecov$ -- covariance
         }

 Compute $\msqrt{R}$ (\ref{equ:sqrtr})\\
 $d \leftarrow \text{measurement dimension}$ \\
 \While{ $\measdim > 0$}
 {
 Compute $\msqrt{\statecov}$ (\ref{equ:sqrtp})\\
 Compute $M$ and $Q_i, 1\leq i \leq \measdim$ (\ref{equ:M}-\ref{equ:Q})\\
 Compute $\trvec$ and $\trmat$  (\ref{equ:Q1}) and (\ref{equ:QQ}) \\
 Compute $\eigVECTOR$ and $\eigVALUE$  (\ref{equ:eigdecom}) \\
 $\decmat \leftarrow \eigVECTOR^T\msqrt{\meascov}^{-1} $ \\
 Choose largest $k$ so that $\eigVALUE\el{i,i}\leq\nonllimit, i \leq k \wedge \eigVALUE\el{j,j}>\nonllimit, j>k$ \\
 \If{$k==0$}{$k\leftarrow 1$}
 	\tcp{Compute partial \ac{EKF2} update}
	\Indp
 	$\predmeas \leftarrow \decmat\el{1:k,:} \left[\measfun(\mean) + \frac{1}{2} \trvec \right] $ \\
	$\innocov \leftarrow \decmat\el{1:k,:}MM^T\decmat\el{1:k,:}^T+ \frac{1}{2} \Lambda\el{1:k,1:k} +  I $	\\
	$\kalmangain \leftarrow \msqrt{\statecov}M^T\decmat\el{1:k,:}^T \innocov^{-1}$ \\
	$\mean \leftarrow \mean + \kalmangain ( \decmat\el{1:k,:}\meas - \predmeas)$ \\
	$\statecov \leftarrow \statecov - \kalmangain \innocov \kalmangain^{T}$ \\
	\Indm
 	\tcp{Update remaining measurement}
	\Indp
	$\meas \leftarrow  \decmat\el{k+1:\measdim,:}\meas$ \\
	$\measfun(\state) \leftarrow  \decmat\el{k+1:\measdim,:} \measfun(\state)$ \\
	$\msqrt{R} \leftarrow I$ \tcp{Updated measurement noise covariance is an identity matrix due to decorrelation}
	$d\leftarrow d-k$ \tcp{Updated measurement dimension}
	\Indm
}
\vspace{0.5ex}
\caption{Algorithm for doing the measurement update in \ac{OEKF2} }
\label{algo:OEKF}
\end{algorithm}
Figure~\ref{fig:pdfs} shows how \ac{OEKF2} treats a two-dimensional second order polynomial measurement function
\begin{equation}
y=\matr{x^2 - 2x - 4 \\
   -x^2 + \frac{3}{2}} + \varepsilon,
\end{equation}
where $\varepsilon \sim \N(0,I)$. The prior has mean $1$ and covariance $1$. The nonlinearity of each measurement  is $4$ and the total nonlinearity is $8$. Then $D=\frac{1}{\sqrt{2}}\matr{1 & 1 \\ 1 &-1}$ and the transformed measurement model has a linear term and a polynomial term
\begin{equation}
\hat y= \sqrt{2} \matr{-x - \frac{5}{4}   \\
   x^2 - x - \frac{11}{4}} + \hat\varepsilon,
\end{equation}
where $\hat\varepsilon \sim \N(0,I)$. After transformation the first element of the measurement function is linear and $\eta_1=0$ and all the nonlinearity is associated with the second element $\eta_2=8$. In \ac{OEKF2} the linear measurement function is applied first and the partially updated state has mean $-\frac{1}{2}$ and covariance $\frac{1}{3}$. The polynomial measurement function is applied using this partially updated state. The amount of nonlinearity for the second order polynomial has decreased from $8$ to $\frac{8}{9}$. \ac{EKF2} applies both measurements at once and the posterior estimate is the same for the original and transformed measurement models as shown in Appendix~\ref{app:lintrans}. When comparing to the true posterior, which is computed using a dense grid, the posterior estimate of \ac{OEKF2} is significantly more accurate than the \ac{EKF2} posterior estimate.

\section{Tests}
\label{sec:tests}

We compare the proposed \ac{OEKF2} with other Kalman filter extensions and a \ac{PF} in three different test scenarios.
The  \ac{OEKF2} was tested with 4 different values for $\nonllimit$. When $\nonllimit=\infty$  the whole measurement is applied at once and the algorithm is a numerical \ac{EKF2}. When $\nonllimit<0$ measurement elements are processed one at a time and when $\nonllimit=0$ all linear measurement elements are first processed together and then nonlinear measurement elements one by one. Due to numerical roundoff errors it is better to use a small positive  $\nonllimit$ to achieve this kind of behaviour. In our tests we use values $\{-\infty, 0.1, 1, \infty\}$ for $\nonllimit$.

\ac{EKF} and \ac{EKF2} are implemented as explained in Section~\ref{sec:numup} with analytical Jacobians and Hessians. \ac{RUF} is implemented according to $\cite{RUKF}$ with 3 and 10 steps. \ac{IKF} uses 10 iterations. For \ac{UKF} the values for sigma point parameters are $\alpha=10^{-3}, \kappa=0, \beta=2$. All  Kalman filter extensions are programmed in Matlab with  similar levels of code optimizations, but the runtimes should still be considered to be only indicative.

For reference we computed estimates with a bootstrap particle filter that does systematic resampling at every time step \cite{systematicresampling} using various numbers of particles and with a \ac{PF} that uses \ac{EKF} for computing the proposal distribution \cite{de2000sequential} with 10 particles.

In every test scenario the state transition model is linear time-invariant $\state_t = \jacobianf \state_{t-1} + \statenoise$, where $\statenoise \sim \N(0,\stateerrcov)$. Thus, the prediction step (\ref{equ:pred1})-(\ref{equ:pred2}) can be computed analytically and all Kalman filter extensions in tests use the analytical prediction.
 
The first test scenario is an aritificial example chosen to show the maximal potential of \ac{OEKF2}. The measurement model used is
\begin{equation}
	\measfun(\state) = \matr{
	2x\el{1} + x\el{2} + x\el{3} + \frac{1}{2}x^2\el{1} + \frac{1}{2}x^2\el{2} + \frac{1}{2}x^2\el{3}\\
	x\el{1} + 2x\el{2} + x\el{3} + \frac{1}{2}x^2\el{1} + \frac{1}{2}x^2\el{2} + \frac{1}{2}x^2\el{3}\\	
	x\el{1} + x\el{2} + 2x\el{3} + \frac{1}{2}x^2\el{1} + \frac{1}{2}x^2\el{2} + \frac{1}{2}x^2\el{3}\\
	x\el{1} + x\el{2} + x\el{3} + x^2\el{1} + \frac{1}{2}x^2\el{2} + \frac{1}{2}x^2\el{3}\\
	x\el{1} + x\el{2} + x\el{3} + \frac{1}{2}x^2\el{1} + x^2\el{2} + \frac{1}{2}x^2\el{3}\\
	x\el{1} + x\el{2} + x\el{3} + \frac{1}{2}x^2\el{1} + \frac{1}{2}x^2\el{2} + x^2\el{3}
		} + \measnoise, \label{equ:measmodel}
\end{equation}
where $\measnoise \sim \N(0, 8I+\mathbf{1})$ and $\mathbf{1}$ is a matrix of ones. This model is a linear transformation of 
\begin{equation}
	\hat \measfun(\state) = \matr{
		x\el{1} \\
		x\el{2} \\
		x\el{3} \\
		\frac{1}{2}x^2\el{1} \\
		\frac{1}{2}x^2\el{2} \\
		\frac{1}{2}x^2\el{3} \\		
		} + \hatmeasnoise, \label{equ:transmeasmodel}
		\end{equation}
		where $\hatmeasnoise \sim \N(0, I)$.  The first three elements of  (\ref{equ:transmeasmodel}) are linear and \ac{OEKF2} with $\nonllimit \in \{0.1,1\}$ uses the three linear measurement functions first to update the state. In this test scenario the prior mean is at the origin, the prior and state transition noise covariances are  both $16I$, and the state transition matrix is an identity matrix.

\begin{figure}
	\includegraphics[width=\columnwidth]{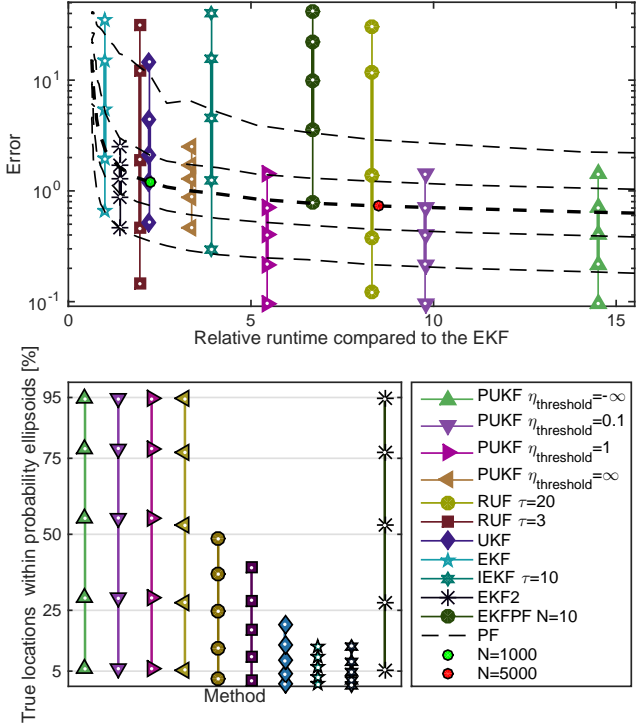}
	\caption{Accuracy of different Kalman filter extensions in estimation with second order polynomial measurement model (\ref{equ:measmodel}). In the top figure markers show the 5\%, 25\%, 50\%, 75\% and 95\% quantiles of errors for each method for every estimated step. The errors are computed as the norm of the difference of the true and estimated mean. In the bottom figure the markers show how often the true state was within estimated error ellipsoids containing 5\%, 25\%, 50\%, 75\% and 95\% of the probability mass. }
	\label{fig:polyres}
\end{figure}

Results for positioning with measurement model (\ref{equ:measmodel}) are presented in Figure~\ref{fig:polyres}. 
 The markers in the upper plot show the 5\%, 25\%, 50\%, 75\% and 95\% quantiles of mean errors for each method. The quantiles are computed  from 1000 runs consisting of 10 steps each. To show the quantiles better a logarithmic scale for error is used. \ac{OEKF2} ($\nonllimit<\infty$) is the most accurate of the Kalman filter extensions by a large margin. When $\nonllimit=\infty$ the whole measurement is processed at once and the result is the same as with \ac{EKF2}, as expected. In this test scenario the \ac{OEKF2} performs clearly the best and methods that use \ac{EKF} linearizations have very large errors. \ac{OEKF2} also outperforms \ac{PF} with similar runtime.

In the bottom plot the accuracy of covariance estimates of different Kalman filter extensions are compared. For this plot we compute how often the true state is within the 5\%, 25\%, 50\%, 75\% and 95\% ellipsoids of the Gaussian posterior. That is, a true location is within the $p$ ellipsoid when
\begin{equation}
	\chi^{2}_\statedim\left((\mean-\truestate)^T\statecov^{-1} (\mean-\truestate) \right) < p,
\end{equation}
where $\truestate$ is the true state, $\mean$ and $\statecov$ are the posterior mean and covariance computed by the filter, and $\chi^{2}_\statedim$ is the cumulative density function of the chi-squared distribution with $\statedim$ degrees of freedom. The filter's error estimate is reliable when markers are close to the $p$ values (dotted lines in the Figure).  From the figure it is evident that \ac{OEKF2} and \ac{EKF2} have the most reliable error estimates and all  other methods have too small covariance matrices. 

The \ac{EKF}\ac{PF} did not perform wery well. This is probably caused by the inconsistency of \ac{EKF} estimates that were used as the proposal distribution. We tested \ac{EKF}\ac{PF} also with 1000 particles. The estimation accuracy was similar to that obtained with a bootrstrap \ac{PF} with 1000 particles, but the algorithm was much slower than other algorithms.

In our second test scenario the planar location of a target is estimated using bearing measurements. When the target is close to the sensor the measurement model is nonlinear, but when the target is far away the measurement becomes almost linear. The measurement model is
\begin{equation}
	\meas = \atan( \state\el{2}-\BSloc\el{2},\state\el{1}-\BSloc\el{1}) + \measnoise, \label{equ:anglemeas}
\end{equation}
where $\atan$ is the four quadrant inverse tangent, $\BSloc$ is the sensor location, and measurement noises are zero mean independent, with standard deviation of $2\degree$. We choose the branch of $\atan$ so that evaluated values are as close as possible to the realized measurement value. In the test scenario two bearings measurements are used, one from a sensor close to the prior and the second from a sensor far away.

\begin{figure}
	\includegraphics[width=\columnwidth,clip=true,trim=1.5cm 1cm 1cm 0cm]{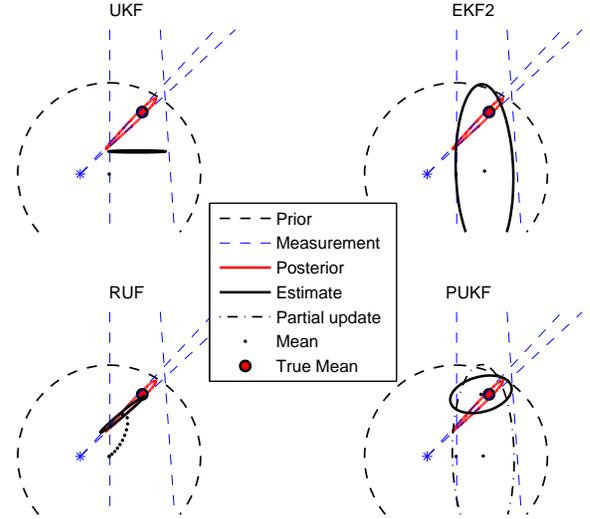}
	\caption{Example situation of bearings positioning}
	\label{fig:bearex}
\end{figure}

A representative initial state update using  \ac{UKF}, \ac{EKF2}, \ac{RUF} and \ac{OEKF2}  is shown in Figure~\ref{fig:bearex}.
This example is chosen so that the differences between estimates of different filters is clearly visible.
The red line encloses the same probability mass of the true posterior as the $1\cdot\sigma$ ellipses (black lines) of the Gaussian approximations computed with different Kalman filter extensions. The measurement from the distant sensor is almost linear within the prior and \ac{UKF} uses it correctly, but the linearization of the estimate from the nearby sensor is not good and the resulting posterior is very narrow (\ac{EKF} would be similar). In the \ac{EKF2} update the second order term of the measurement model from the nearby sensor is so large that \ac{EKF2} almost completely ignores that measurement and the prior is updated using only the measurement from the distant sensor. The iterative update of \ac{RUF} results in an estimate with small covariance that has similar shape as the true covariance. The mean of the true posterior is not inside the one-sigma ellipses of the $\ac{RUF}$ estimate and the mean is too close to the nearby sensor.

The first transformed measurement used by \ac{OEKF2} is almost the same as the measurement from the distant sensor and the estimate after the first partial update is similar to the \ac{EKF2} estimate. Because the estimate updated with the first measurement is further away from the nearby sensor the linearization of the second measurement is better and the posterior estimate is closer to the true posterior than with \ac{EKF2}. The covariance estimate produced by \ac{OEKF2} is more conservative than the \ac{RUF} of \ac{UKF} covariances. 

\begin{figure}
	\includegraphics[width=\columnwidth]{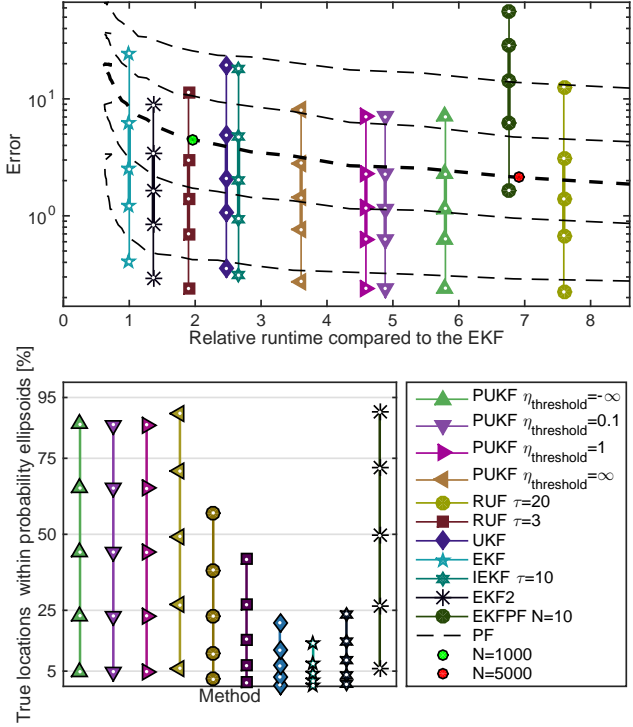}
	\caption{Accuracy of different filters in bearings only tracking}
	\label{fig:bearres}
\end{figure}

Figure~\ref{fig:bearres} shows the statistics for this scenario. For this Figure the scenario was ran 1000 times using the same sensor locations and 10 step estimation with a 4-dimensional state model containing 2 position and 2 velocity dimensions.  The  prior has zero mean and covariance $10I$. 
The state transition function is 
\begin{equation}
\statefun(x)= \matr{I & I \\ \mathbf{0} & I } x+\statenoise, \label{equ:modold}
\end{equation}
where 
\begin{equation}
\statenoise \sim
\N\left(0, \matr{\frac{1}{300} I  & \frac{1}{200} I\\ \frac{1}{200}I & \frac{1}{100}I}\right).
\end{equation}

 Figure~\ref{fig:bearres} shows that the \ac{OEKF2} provides the best accuracy. Interestingly \ac{RUF} with 3 iterations has better accuracy than with $20$ iterations. From the plot that shows the accuracy of the error estimates we can see that the \ac{OEKF2} and \ac{EKF2} have the best error estimates. Other methods have too optimistic covariance estimates. In this test scenario the \ac{PF} did not manage to get good estimates with similar runtimes. 

In the third test scenario we consider bearings only tracking with sensors close to each other. Otherwise the measurement model is the same as in the previous scenario
The prior is  as in previous test scenario. 
The state transition function is also (\ref{equ:modold})
but the state transition noise is higher:
\begin{equation}
\statenoise \sim
\N\left(0, \matr{\frac{1}{3} I  & \frac{1}{2} I\\ \frac{1}{2}I & I}\right).
\end{equation}

\begin{figure}
	\includegraphics[width=\columnwidth,clip=true,trim=1.5cm 1cm 1cm 0cm]{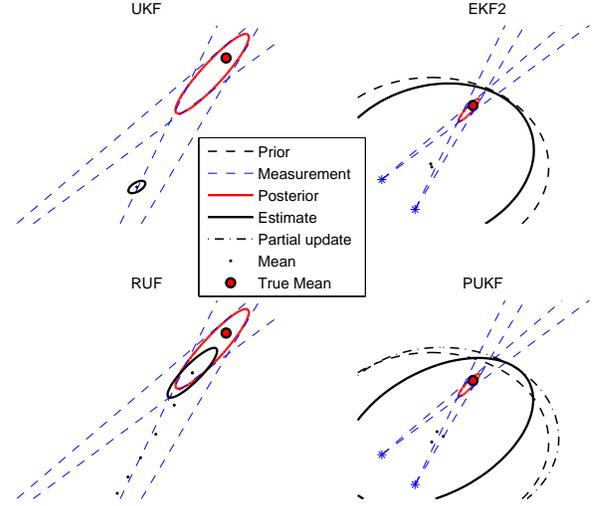}
	\caption{Example of first update in bearings only tracking}
	\label{fig:bear2}
\end{figure}

The initial state and representative first updates are shown in Figure~\ref{fig:bear2}.
In this Figure \ac{UKF} and \ac{RUF} estimates have very small covariances and so the plots are magnified. The \ac{UKF} estimate mean is closer to the true mean than \ac{EKF2} and \ac{OEKF2} estimates, but the  covariance of the estimate is very small. \ac{RUF} has a better estimate than \ac{UKF}, but the estimate is biased towards the sensor locations. Because both sensors are nearby and have large second order terms \ac{EKF2} and \ac{OEKF2} estimates do not differ much.

\begin{figure}
	\includegraphics[width=\columnwidth,clip=true,trim=0cm 0cm 0cm 0cm]{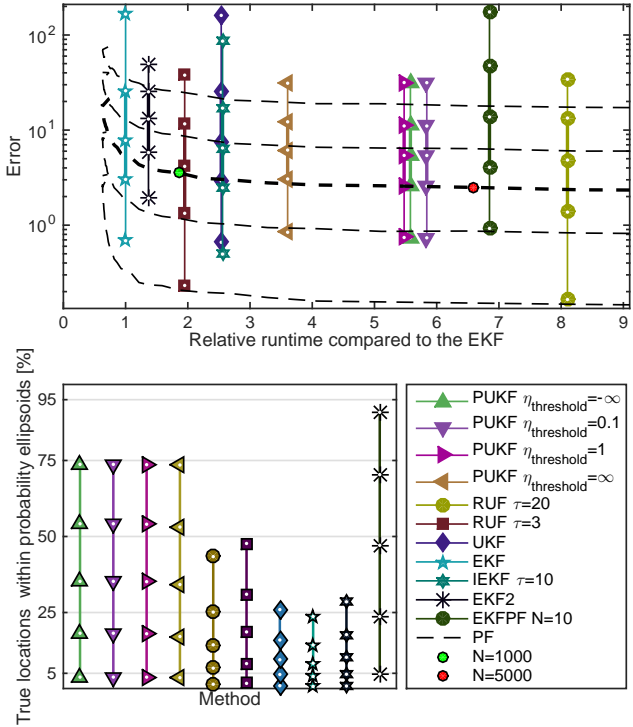}
	\caption{Results for bearings only tracking with sensors close to each other}
	\label{fig:bearres2}
\end{figure}

Results for estimating 10 step tracks 1000 times are shown in Figure~\ref{fig:bearres2}. In this case the \ac{RUF} has the best accuracy. In \ac{OEKF2} there is only very small differences whether all of the measurement are used at once or a nonlinearity threshold is used. This means that in this measurement geometry the partitioned update does not improve  accuracy. \ac{EKF2} has better covariance estimates than the numerical update \ac{OEKF2} even though it has larger errors. The covariance estimates produced by \ac{RUF} were again too small. In this test the \ac{PF} has better accuracy than the Kalman filter extensions. 

To further evaluate the accuracy of the estimates, we compare \ac{KL} divergences of estimates. The \ac{KL} divergence is defined as
\begin{equation}
	\int \ln\left( \frac{p(\state)}{q(\state)}  \right) p(\state) \mathrm{d}\state. \label{equ:KL},
\end{equation}
where $p(x)$ is the \ac{pdf} of the true distribution and $q(x)$ is the \ac{pdf} of the approximate distribution \cite{Kullback51klDivergence}.
We computed the \ac{KL} divergence for position dimensions. The true \ac{pdf} is approximated using a $50\times50$ grid. The probability for each grid is computed as the sum of particle weights of a \ac{PF} particles within each cell. For this we used $10^6$ particles.
\begin{table*}
\centering
\caption{Median Kullback-Leibler divergences of position dimensions in the two bearings only tests}
\label{tbl:KLdiv}
\begin{tabular}{r|c|c|c|c|c|c|c|c|c|c}
Method & PUKF & PUKF & PUKF & PUKF & RUF & RUF  & UKF & EKF& IEKF & EKF2  \\ 
 &  $\nonllimit=-\infty$ &  $\nonllimit=0.1$&  $\nonllimit=1$ &  $\nonllimit=\infty$ &  $\tau=20$ &  $\tau=3$ &  & &  $\tau=10$ &   \\ \hline
First bearings test&0.62 & 0.62 & 0.63 &1.07 & 0.99 &1.70 & 5.20& 14.42&  6.69 & 1.16 \\
Second bearings test &2.14 & 2.14 & 2.14 & 2.33 & 3.56 & 2.53 & 9.36 & 10.97&  8.84 & 2.60
\end{tabular}
\end{table*}

Table~\ref{tbl:KLdiv} shows the median Kullback-Leibler divergences for each method in the two bearings measurement test scenarios. \ac{OEKF2} has the smallest \ac{KL} divergence in both test scenarios.

\section{Conclusions and Future Work}
\label{sec:conclusions}
In this paper we presented a new extension of the Kalman filter: the \acf{OEKF2}. The  proposed filter evaluates the nonlinearity of a multidimensional measurement and transforms the measurement model so that some dimensions of the measurement model have as low nonlinearity as possible. \ac{OEKF2} does the update of the state using the measurement in parts, so that the parts with the smallest amounts of nonlinearity are processed first. The proposed algorithm improves estimation results when measurements are such that the partial update reduces the nonlinearity of the remaining part. According to the simulated tests the \ac{OEKF2} improves the estimates when measurements can be transformed so that an informative linear part of the measurement can be extracted.

In many practical situations the almost linear part could be extracted manually. For example, \ac{GPS} measurements are almost linear and they could be applied before other measurements. The proposed algorithm does the separation automatically and when using the numerical algorithm for computing the prediction and update analytical differentiation is not required.

In our tests the estimated covariances produced by \ac{EKF2} and \ac{OEKF2} were the most accurate. In \cite{RUKF} it was claimed that \ac{RUF} produces more accurate error estimates than \ac{EKF2}. Their results were based on comparing $3\cdot\sigma$ errors in 1D estimation. In this comparison 92\% of samples should be within the $3\cdot\sigma$ range. For their results they had only 100 samples and from the resulting figure it is hard to see how many samples exactly are within the range, but for \ac{EKF2} most of the points are within the range and some are outside. 

In our tests, among other Kalman filter extensions \ac{RUF} had good accuracy, but it provided too small covariance matrices. In future it could be interesting to extend \ac{RUF} \cite{RUKF} to use \ac{EKF2}-like statistical second order linearization and then combine it with the proposed algorithm.

Another use case for \ac{OEKF2} would be merging it with the Binomial Gaussian mixture filter \cite{binoGMF}. This filter decorrelates measurements and uses nonlinearity measure (\ref{equ:ithnonl}) as an indication of whether the measurement model is so nonlinear that the prior component should be split. By decorrelating measurements with the algorithm proposed in this paper and doing the partial updates for the most linear components first, unnecessary splits could be avoided.

\bibliographystyle{myIEEEtran}
\bibliography{viitteet_pukf}

\appendices

\section{Invariance of \ac{EKF} and \ac{EKF2} to a linear transformation of the measurement model}
\label{app:lintrans}

The second order Taylor polynomial approximation of the measurement function is
\begin{equation}
\begin{aligned}
	 \measfun(\state) = & \measfun(\priormean) + \jacobian (\state - \priormean)  \\ & + \frac{1}{2} \matr{ (\state - \priormean)^T \hessian_1 (\state - \priormean) \\ (\state - \priormean)^T \hessian_2 (\state - \priormean) \\ \vdots \\ (\state - \priormean)^T \hessian_\measdim (\state - \priormean)}  + \measnoise \label{equ:secondorder}
	 \end{aligned}
\end{equation}
where Jacobian $\jacobian$ and Hessians $\hessian$ are evaluated at prior mean, $\measnoise$ is the measurement function noise.

In the linear transformation the measurement function (\ref{equ:secondorder}) is multiplied by $\decmat$. The second order approximation is
\begin{equation}
\begin{aligned}
\hat{\measfun}(\state) = & \decmat\measfun(\state) = \decmat\measfun(\priormean) + \decmat\jacobian(\state-\priormean) \\ &+ \frac{1}{2}\decmat \matr{(\state-\priormean)^T \hessian_1 (\state-\priormean) \\ (\state-\priormean)^T \hessian_2 (\state-\priormean) \\ \vdots \\ (\state-\priormean)^T \hessian_\statedim (\state-\priormean)  } + \decmat\measnoise 
\end{aligned}
\end{equation}
The transformed Jacobian is 
\begin{align}
\hat{\jacobian}&= \decmat\jacobian 
\end{align}
and $i$th transformed  Hessian is  
\begin{align}
\hathessian_i&= \sum_{k=1}^\statedim \decmat\el{i,k}\hessian_k.
\end{align}
The terms $\trvec$ and $\trmat$ are 
\begin{align}
&\begin{aligned}
\hattrvec &= \matr{\tr \priorcov {\hathessian}_1 \\ \tr \priorcov {\hat\hessian}_2 \\ \vdots \\ \tr \priorcov {\hat\hessian}_\statedim}  = \matr{ \tr \priorcov \sum_{k=1}^\statedim \decmat\el{1,k} \hessian_k \\  \tr \priorcov \sum_{k=1}^\statedim \decmat\el{2,k} \hessian_k \\ \vdots \\  \tr \priorcov \sum_{k=1}^\statedim \decmat\el{\statedim,k} \hessian_k} \\
&= \decmat \matr{\tr \priorcov\hessian_1\\ \tr \priorcov\hessian_2\\ \vdots \\ \tr \priorcov\hessian_\statedim }= \decmat \trvec
\end{aligned}
\\
&\begin{aligned}{\hattrmat}\el{i,j} &= \tr \priorcov {\hathessian}_i \priorcov {\hathessian}_j  \\
& =  \tr \priorcov \left(  \sum_{k=1}^\statedim \decmat\el{i,k}\hessian_k \right) \priorcov \left(  \sum_{l=1}^\statedim \decmat\el{i,l}\hessian_l \right)\\
&=    \sum_{k=1}^\statedim  \sum_{l=1}^\statedim \decmat\el{i,k}   \decmat\el{j,l} \tr \priorcov \hessian_k \priorcov\hessian_l \\
 \Rightarrow {\hattrmat} & = \decmat\trmat \decmat^T \\
\end{aligned}
\end{align}
For \ac{EKF} update these terms are replaced with zero matrices.

Now using these transformed quantities in the \ac{EKF2} update equations (\ref{equ:predmeas}-\ref{equ:ekfend}) gives
\begin{align}
&\begin{aligned}
	\hatpredmeas & = \hat \measfun(\priormean)+ \frac{1}{2} \hattrvec &= \decmat\left(\measfun(\priormean)+ \frac{1}{2} \trvec \right)  \\
\end{aligned} \\
&\begin{aligned}
\hat{\innocov}&=\decmat\jacobian\priorcov\jacobian^T\decmat^{T}+\frac{1}{2}\decmat\trmat \decmat^T+ \decmat \meascov \decmat^T \\
&=\decmat\innocov \decmat^T 
\end{aligned} \\
&\begin{aligned}
\hat{\kalmangain}&=\priorcov \hat{\jacobian}^T\hat{\innocov}^{-1} \\
& = \priorcov \jacobian^T\decmat^{T}\decmat^{-T}\innocov^{-1}\decmat^{-1}   =\priorcov \jacobian^T\innocov^{-1}\decmat^{-1} \\
& = \kalmangain \decmat^{-1}
	 \end{aligned}\\
&\begin{aligned}
\hat{\mean}^+&=\priormean + \hat{\kalmangain}\left[\decmat\meas-\decmat \measfun(\priormean) - \frac{1}{2} \decmat \trvec \right] \\& =\priormean + \kalmangain(\meas-\jacobian\priormean - \frac{1}{2}\trvec) \\ & = \mean^+   
	 \end{aligned}\\
&\begin{aligned}
\hat{\statecov}^+&=\priorcov-\hat{\kalmangain}\hat{\innocov}\hat{\kalmangain}^T \\ &  = \priorcov - \kalmangain \decmat^{-1}\decmat\innocov \decmat^T (\kalmangain \decmat^{-1})^T   =\priorcov - \kalmangain\innocov\kalmangain^T  \\
& = \statecov^+,
\end{aligned}
\end{align}
which shows that the posterior is the same as with the non-transformed measurements.
\section{Proof that the nonlinearities are minimized}
\label{app:nonlminimization}
Let $\trmat$ be a diagonal matrix containing nonlinearity values ordered ascending on the diagonal and let measurement noise covariance matrix be identity matrix $\meascov=I$. We will show that the smallest diagonal element of $\trmat$ is as small as possible under a linear transformation that preserves  $\meascov=I$ and further that the second smallest diagonal element is as small as possible, when the smallest is as small as possible etc.
\balance

If the measurement model is transformed by multiplying it with $V$, the transformed variables are $\hat \trmat= V \trmat V^T$ and $\meascov=VIV^T=VV^T$. Because we want to have $\meascov=I$, $V$ has to be unitary. The $i$th diagonal element of the transformed matrix is $v_i^T \trmat  v_i = \sum_{j=1}^\measdim v_{i,[j]}^2 \trmat\el{j,j}$, where $v_i$ is the $i$th column of $V$. Because $V$ is unitary, we have $\sum_{j=1}^\measdim v_{i,[j]}^2=1$ and the $i$th diagonal element of the transformed matrix $\hat\trmat$ is
\begin{equation}
\sum_{j=1}^d v_{i, [j]}^2 \trmat_{[j,j]} \geq \sum_{j=1}^d v_{i, [j]}^2 \min_j\{\trmat_{[j,j]}\} = \min_j \{\trmat_{[j,j]}\}.
\end{equation}
Thus, the new diagonal element cannot be smaller than the smallest diagonal element of $\trmat$.

If the smallest element is in the first element of the diagonal the possible transformation for the second smallest element is
\begin{equation}
	\hat{\trmat} = \matr{1 & 0^T \\ 0  & V} \trmat \matr{1 & 0^T \\ 0  & V^T} .
\end{equation}
With the same reasoning as given already the second diagonal has to be already the smallest possible. Inductively this applies to all diagonal elements.

\end{document}